\documentclass[11pt,a4paper]{article}
\usepackage{authblk}
\usepackage[utf8]{inputenc}
\usepackage[T1]{fontenc}
\usepackage{graphicx}
\usepackage{float}
\usepackage{subfigure}
\usepackage{multirow}
\usepackage{amssymb}
\usepackage{amsmath}
\usepackage{amsfonts}
\usepackage[numbers,sort&compress]{natbib}
\usepackage{hyperref}
\usepackage{caption}
\usepackage{subfigure}
\usepackage{mathrsfs}
\usepackage[title]{appendix}
\usepackage{xcolor}
\usepackage{textcomp}
\usepackage{manyfoot}
\usepackage{booktabs}
\usepackage{algorithm}
\usepackage{algorithmicx}
\usepackage{algpseudocode}
\usepackage{listings}
\allowdisplaybreaks
\newtheorem{thm}{Theorem}[section]
\newtheorem{deff}{Definition}[section]

\newtheorem{ex}{Example}[section]
\numberwithin{deff}{section}
\numberwithin{thm}{section}

\title{Analysis of the maps with variable fractional order}
\author[1]{Prashant M. Gade}
\author[2]{Sachin Bhalekar\footnote{Corresponding author email sachinbhalekar@uohyd.ac.in}}
\author[2]{Janardhan Chevala}
\affil[1]{Department of Physics, Rashtrasant Tukadoji Maharaj Nagpur University, Nagpur, India.}
\affil[2]{School of Mathematics and Statistics, University of Hyderabad, Hyderabad, 500046 India.}

\date{}
\begin{document}
	
	\maketitle

		\begin{abstract}
	Fractional order differential and difference  equations are used to model systems with memory. Variable order fractional equations are proposed to model systems where the memory changes in time. We investigate stability conditions for linear variable order difference equations where the order is periodic function with period $T$. We give a general procedure for arbitrary $T$ and for $T=2$ and $T=3$, we give exact results. For $T=2$, we find that the lower order determines the stability of the equations. For odd $T$, numerical simulations indicate that we can approximately determine the stability of equations from the mean value of the variables.	
	\end{abstract}
	
	
\section{Introduction} 
Memory and non-locality of the system can change with time. Such non-locality can be modeled by non-stationary power-law kernel \cite{lorenzo2002variable,sun2019review}. 
In 1993, Samko and Ross \cite{samko1993integration} first introduced the idea of variable order differential and integral, along with a few fundamental characteristics. Derakhshan \cite{derakhshan2022existence} considered the Caputo type variable order fractional differential equation (VO-FDE) in fluid mechanics and discussed Ulam–Hyers stability conditions and  the existence–uniqueness of a solution of the Caputo type VO-FDE. Sarvar \cite{sarwar2022existence} considered Caputo type variable order fractional initial value problem, established the global existence of solution and discussed different types of Ulam–Hyers stability results. Albasheir et al. \cite{albasheir2023theoretical} proposed Ulam–Hyers stability conditions and existence–uniqueness for the Caputo FDEs of variable order. Such results are under strict conditions. The functions satisfying all their conditions are rare, hence the results are of limited use. Lyapunov stability (of equilibrium points) of variable order fractional differential equations (VOFDE) is not available in the literature, to the best of our knowledge.

\par The numerical solutions for VO-FDE can be found in \cite{solis2018novel,zeng2015generalized,zhao2015second,shen2012numerical,chen2014numerical,cao2016high,wu2022fractional}. The initial value problem of the nonlinear Riemann-Liouville VO-FDE is numerically solved by Wang et al. \cite{wang2019analysis} utilizing a finite difference approach. An exponentially accurate fractional spectral collocation approach for solving linear/nonlinear VO-FDEs is developed by Zayernouri and Karniadakis \cite{zayernouri2015fractional}. The numerical solution of a Caputo-type VO-FDE system based on shifted Legendre polynomials is examined by Derakhshan \cite{derakhshan2022existence}. The variable order fractional differential equations are widely used in the fields Biology \cite{khan2023system}, Mechanics \cite{coimbra2003mechanics}, Finance \cite{soradi2020king} etc. The detailed review of the progress in this field is provided by Sun et al. in \cite{sun2019review} and Patnaik et al. in \cite{patnaik2020applications}.    

\par On the other hand, the fractional variable-order discrete maps are currently of interest to researchers. Mozyrska et al. \cite{mozyrska2019stability} considered the Grünwald-Letnikov type linear fractional variable order discrete-time system and discussed the stability conditions. Hamadneh et al. \cite{hamadneh2023finite} considered the Caputo type non-linear nabla variable-order fractional discrete neural network and discussed finite time stability conditions. Wu et al. \cite{wu2019new} investigated a new variable-order fractional chaotic systems and applied to block image encryption. There has been work on stability analysis of Grunwald-Letnikov type linear variable order fractional equations \cite{mozyrska2019stability}. The stability conditions are different from those derived in this work because we study Caputo type fractional difference equations.
Such systems have been numerically investigated \cite{bensid2022discrete,almatroud2024new}. To the best of our knowledge, the stability analysis of variable-order Caputo type fractional difference equations is not done by any researcher. This motivated us to present the stability results on such maps with periodic-order.  We present exact results for the periodic order of period-2 and period-3. We follow with a numerical exploration
of higher periods. We find that the mean value of orders plays an important role for odd values of higher periods.
 
\par The paper is organized as follows:
In Section \ref{prel}, we provide a few key definitions and theorems. In Sections \ref{per2} and \ref{per3}, we discuss the stability region of considered system with periodic order two and three, respectively, and provide some examples. In Section \ref{pern}, we present a numerical exploration of higher periods of considered system. Finally, the results are summarized in Section \ref{con.}.

\section{Preliminaries} \label{prel}
In this section, we present some basic definitions and results.\\ 
Let $h>0$, \,$ a\in \mathbb{R}$,
$(h\mathbb{N})_a=\{a,a+h,a+2h,\dots\}$ and $\mathbb{N}_a=\{a,a+1,a+2,\dots\}$.
\begin{deff}(See\cite{bastos2011discrete,ferreira2011fractional,mozyrska2015transform})
	For a function $x:(h\mathbb{N})_a\longrightarrow \mathbb{R}$, the forward h-difference operator is defined as
	$$(\Delta_hx)(t)=\frac{x(t+h)-x(t)}{h},$$
	where $t\in(h\mathbb{N})_a$. \\
	Throughout this paper, we take $a=0$ and $h=1$. 
\end{deff}

\begin{deff}\cite{mozyrska2015transform}
	For a function $x:\mathbb{N}_0\longrightarrow \mathbb{R}$, the fractional sum of order $\beta >0$\, is given by 
	$$(\Delta^{-\beta} x)(t)=\frac{1}{\Gamma(\beta)} \sum_{s=0}^n{\frac{\Gamma(\beta+n-s)}{\Gamma(n-s+1)}}x(s),$$  
	where $t=\beta +n$, $n\in \mathbb{N}_0$.
\end{deff}

\begin{deff}\cite{fulai2011existence,mozyrska2015transform}  
	Let $\mu>0$ and $m-1<\mu<m$, where $m\in\mathbb{N}.$
	The $\mu$th fractional Caputo-like difference is defined as
	$$\Delta^\mu x(t)=\Delta^{-(m-\mu)}(\Delta^mx(t)),$$
	where $t\in\mathbb{N}_{m-\mu}$ and
	$$\Delta^mx(t)=\sum_{k=0}^m \left(\begin{array}{c}m\\k\end{array}\right)(-1)^{m-k}x(t+k).$$
\end{deff}

\begin{deff}\cite{mozyrska2015transform} 
	The Z-transform of a sequence $\{y(n)\}_{n=0}^\infty$ is a complex function given by $Y(z)=Z[y](z)=\sum_{k=0}^\infty y(k)z^{-k}$, where $z\in \mathbb{C}$ is a complex number for which the series converges absolutely.
\end{deff}
\begin{deff}\cite{mozyrska2015transform} 
	Let 
    \begin{eqnarray*}
		\Tilde{\phi}_\alpha(n) &=& \frac{\Gamma(n+\alpha)}{\Gamma(\alpha)\Gamma(n+1)} \nonumber \\
		&=& \left(\begin{array}{c} n+\alpha-1\\n \end{array} \right) = (-1)^n \left( \begin{array}{c}
			-\alpha \\ n \end{array} \right)
	\end{eqnarray*}
    be a family of binomial functions defined on $\mathbb{Z}$, parametrized by $\alpha$.\\
Then, 
	\begin{eqnarray*}
		Z(\Tilde{\phi}_\alpha(t))=\frac{1}{(1-z^{-1})^\alpha},\; |z|>1. \nonumber
	\end{eqnarray*}
\end{deff}
\begin{deff}\cite{mozyrska2015transform}
	The convolution $\phi\ast x$ of the functions $\phi$ and $x$ defined on $\mathbb{N}$ is defined as 
	\begin{eqnarray*}(\phi\ast x)(n)=\sum_{s=0}^n \phi(n-s)x(s)=\sum_{s=0}^n \phi(s) x(n-s).
	\end{eqnarray*} 
	Then, the Z-transform of this convolution is
	\begin{eqnarray*}
		Z(\phi\ast x)(n)=(Z(\phi)(n)) Z((x)(n)).
	\end{eqnarray*} 
\end{deff}
\begin{thm}\cite{fulai2011existence}
	The difference equation $$\Delta^{\alpha} x(t)=f(x(t+\alpha-1))-x(t+\alpha-1),$$ where $0<\alpha \leq 1$, $t \in \mathbb{N}_{1-\alpha}$, 
	is  equivalent to
	\begin{equation} 
		x(t)=x(0)+\frac{1}{\Gamma(\alpha)} \sum_{j=1}^t \frac{\Gamma(t-j+\alpha)}{\Gamma(t-j+1 )} (f(x(j-1))-x(j-1)), \label{1}
	\end{equation}
	where $t \in \mathbb{N}_{0}$.
\end{thm}

\begin{deff}\cite{elaydi2005systems,hirsch2012differential}
	A steady state solution or an equilibrium $x_*$ of (\ref{1}) is a real number satisfying $f(x_*) = x_*$.
\end{deff}

\begin{deff} \cite{elaydi2005systems,hirsch2012differential}
	An equilibrium $x_*$ is stable if for each  $\epsilon>0$, there exists $\delta>0$ such that $|x_0 -x_* | < \delta $ implies
	$|x(t) - x_* | < \epsilon$, $t=1,2,3,...$\\
	If $x_*$ is not stable then it is unstable.
\end{deff}

\begin{deff} \cite{elaydi2005systems,hirsch2012differential}
	An equilibrium point $x_*$  is asymptotically stable  if it is stable and there exists $\delta>0$ such that $|x_0 -x_* | < \delta $ implies $ lim_{t\to\infty}x(t) =x_*$.
\end{deff}
Consider the linear map as 
\begin{equation}
	X(t+1)=X(0)+\sum_{j=0}^t \frac{\Gamma(t-j+\alpha)}{\Gamma(\alpha)\Gamma(t-j+1)}(A-I)X(j), \label{2}
\end{equation}
where A is $N \times N$  real matrix, I is $N \times N$ identity matrix and $X(t)\in \mathbb{R}^N $ for each t=0,1,2,...
\begin{thm} \cite{fulai2011existence,abu2013asymptotic,vcermak2015explicit,stanislawski2013stability,bhalekar2022stability}.
	The zero solution of the system (\ref{2}) is asymptotically stable if and only if all the roots of the characteristic equation  $det(z(1-z^{-1})^\alpha I-(A-I))=0$ satisfy $|z|<1$.
\end{thm}
In this work, we focus on the map (\ref{1}) with variable order $\alpha(t)$. We consider
\begin{equation} 
	x(t+1)=x(0)+\frac{1}{\Gamma(\alpha(t))} \sum_{j=0}^t \frac{\Gamma(t-j+\alpha(t))}{\Gamma(t-j+1 )} (r -1)x(j), \label{3}
\end{equation}
where $r \in \mathbb{R}$.
\section{Maps with period-2 fractional order} \label{per2}
Consider the map (\ref{3}) with 
\begin{equation*}
	\alpha(t)=
	\begin{cases}
		\alpha_1,\quad  \text{if t is odd integer} \\
		\alpha_2, \quad \text{if t is even integer} 
	\end{cases}
\end{equation*}
where $0 < \alpha_1 < \alpha_2 \le 1$, $t \in \mathbb{N} \cup \{0\}$.\\
Thus, $\alpha$ is a periodic function of $t$.\\
Let us define $a(t)=x(2t)$ and $b(t)=x(2t+1)$. \\
The system (\ref{3}) can be written in an equivalent form as 
\begin{equation}
	\begin{split}
		a(t+1) &= x(0)+\sum_{k=0}^t \frac{\Gamma(2t+1-2k+\alpha_1)}{\Gamma(\alpha_1) \Gamma(2t-2k+2)} (r-1) a(k)\\
		& + \sum_{k=0}^t \frac{\Gamma(2t-2k+\alpha_1)}{\Gamma(\alpha_1) \Gamma(2t-2k+1)} (r-1) b(k),\\
		b(t+1) &= x(0)+\sum_{k=0}^{t+1} \frac{\Gamma(2t+2-2k+\alpha_2)}{\Gamma(\alpha_2) \Gamma(2t-2k+3)} (r-1) a(k)\\
		& + \sum_{k=0}^t \frac{\Gamma(2t-2k+1+\alpha_2)}{\Gamma(\alpha_2) \Gamma(2t-2k+2)} (r-1) b(k). \label{4}
	\end{split}
\end{equation}
If we define $\phi_{\alpha_1}(n_0)=\Tilde{\phi}_{\alpha_1}(2n_0+1)$, $\phi_{\alpha_1}^{'} (n_0)=\Tilde{\phi}_{\alpha_1}(2n_0)$, $\phi_{\alpha_2}(n_0)=\Tilde{\phi}_{\alpha_2}(2n_0+2)$ and $\phi_{\alpha_2}^{'} (n_0)=\Tilde{\phi}_{\alpha_2}(2n_0+1)$, then the system (\ref{4}) can be written as 
\begin{equation}
	\begin{split}
		a(t+1) &= x(0)+ (r-1) ((\phi_{\alpha_1}\ast a)(t))+(r-1) ((\phi_{\alpha_1}^{'}\ast b)(t)),\\
		b(t+1) &= x(0)+(r-1)( (\phi_{\alpha_2}\ast a)(t))+(r-1)a(t+1)\\
        &+(r-1)( (\phi_{\alpha_2}^{'}\ast b)(t)). \label{5}
	\end{split}
\end{equation}

We have the following Z-transforms, 
\begin{equation}
	\begin{split}
		Z(\phi_{\alpha_1}(t)) &= \frac{\sqrt{z}[(1-z^{-1/2})^{-\alpha_1}-(1+z^{-1/2})^{-\alpha_1}]}{2},\\
		Z(\phi_{\alpha_1}^{'}(t)) &= \frac{(1+z^{-1/2})^{-\alpha_1}+(1-z^{-1/2})^{-\alpha_1}]}{2},\\
		Z(\phi_{\alpha_2}(t)) &= -z+\frac{z[(1+z^{-1/2})^{-\alpha_2}+(1-z^{-1/2})^{-\alpha_2}]}{2},\\
		Z(\phi_{\alpha_2}^{'}(t)) &= \frac{\sqrt{z}[(1-z^{-1/2})^{-\alpha_2}-(1+z^{-1/2})^{-\alpha_2}]}{2}, \\
		Z(a(t+1)) &= zA(z)-za(0), \, \text{and}\, Z(b(t+1))=zB(z)-zb(0),  \label{6}
	\end{split}
\end{equation}
where $Z(a(t))=A(z)$, $Z(b(t))=B(z)$, $a(0)=x(0)$, $b(0)=x(1)$.\\
Applying Z-transform to the system (\ref{5}) and using (\ref{6}), we get
\begin{equation*}
	\begin{split}
   &(z-(r-1) \frac{\sqrt{z}[(1-z^{-1/2})^{-\alpha_1}-(1+z^{-1/2})^{-\alpha_1}]}{2})A(z)\\
		&-((r-1)\frac{[(1+z^{-1/2})^{-\alpha_1}+(1-z^{-1/2})^{-\alpha_1}]}{2})B(z)=\frac{x_0}{1-z^{-1}}+zx(0),
        \end{split}
\end{equation*}
\begin{equation}
	\begin{split}
    &((r-1)\frac{z[(1+z^{-1/2})^{-\alpha_2}+(1-z^{-1/2})^{-\alpha_2}]}{2})A(z)\\
		&+((r-1)\frac{\sqrt{z}[(1-z^{-1/2})^{-\alpha_2}-(1+z^{-1/2})^{-\alpha_2}]}{2}-z)B(z)\\
        &=(r-1)zx(0)-zx(1)-\frac{x_0}{1-z^{-1}}.
         \label{7}
	\end{split}
\end{equation}		
The characteristic equation of the system (\ref{4}) can be obtained by equating the determinant of coefficients of the terms A(z) and B(z) in the system (\ref{7}) to zero as below :		
\begin{equation}
	\begin{split}
		&\frac{-(r-1)^2[(\sqrt{z}-1)^{-\alpha_2}(\sqrt{z}+1)^{-\alpha_1}+(\sqrt{z}-1)^{-\alpha_1}(\sqrt{z}+1)^{-\alpha_2}]}{(\sqrt{z})^{-\alpha_1-\alpha_2}}\\
		& +\sqrt{z}(r-1) \left[ \frac{(\sqrt{z}+1)^{-\alpha_1}-(\sqrt{z}-1)^{-\alpha_1}}{(\sqrt{z})^{-\alpha_1}}+\frac{(\sqrt{z}+1)^{-\alpha_2}-(\sqrt{z}-1)^{-\alpha_2}}{(\sqrt{z})^{-\alpha_2}}\right]\\
        &+2z=0.		\label{8}
	\end{split}
\end{equation}

\subsection{The Stable Region} \label{3.1}
Each component of the system (\ref{4}) is locally asymptotically stable if and only if all the roots $z$ of the characteristic equation (\ref{8}) satisfy $|z|<1$. Therefore, the boundary of the stable region of the system (\ref{4}) can be obtained by substituting $z=\pm 1$ in the characteristic equation (\ref{8}).\\
Multiplying the characteristic equation (\ref{8}) by $(\sqrt{z}-1)^{\alpha_2}$, we get\\

\begin{equation}
	\begin{split}
		&-(r-1)^2 (\sqrt{z})^{\alpha_1+\alpha_2}[(\sqrt{z}+1)^{-\alpha_1}+(\sqrt{z}-1)^{\alpha_2-\alpha_1}(\sqrt{z}+1)^{-\alpha_2}]\\
		&+\sqrt{z}(r-1)[(\sqrt{z})^{\alpha_1}(\sqrt{z}+1)^{-\alpha_1}(\sqrt{z}-1)^{\alpha_2}-(\sqrt{z})^{\alpha_1}(\sqrt{z}-1)^{\alpha_2-\alpha_1}\\
		&+(\sqrt{z})^{\alpha_2}(\sqrt{z}+1)^{-\alpha_2}(\sqrt{z}-1)^{\alpha_2}-(\sqrt{z})^{\alpha_2}]+2z (\sqrt{z}-1)^{\alpha_2}=0. \label{9}
	\end{split}
\end{equation}
Substituting $z=\sqrt{z}=1$ in equation (\ref{9}), we get\\
\begin{equation}
-(r-1)^2( 2^{-\alpha_1})+(r-1)(-1)=0 \label{10}
\end{equation}
Solving (\ref{10}) for $r$, we get $r=1-2^{\alpha_1}$ as the left boundary and $r=1$ as the right boundary of the stable region of the system (\ref{4}), where $\alpha_1=$ min\{$\alpha_1,\alpha_2$\}.\\
Therefore, the stable region of the system (\ref{4}) is $1-2^{\text{min}\{\alpha_1,\alpha_2\}}<r<1$.

\subsection{Corroboration of Analytic Results using Numerical Simulations} \label{3.2}
\begin{figure}[H]
	\centering
	\subfigure[$r=0.4$.]
	{\includegraphics[height=2.2in,keepaspectratio,width=2.2in]{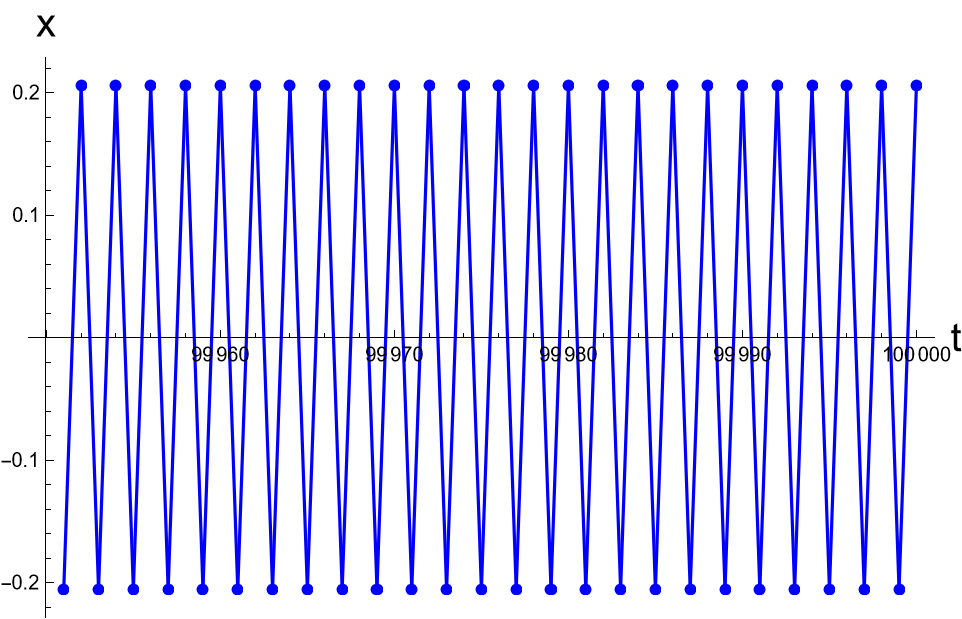}\label{fig1a}} \hspace{0.2cm}
	\subfigure[$r=1.1$.]
	{\includegraphics[height=2.2in,keepaspectratio,width=2.2in]{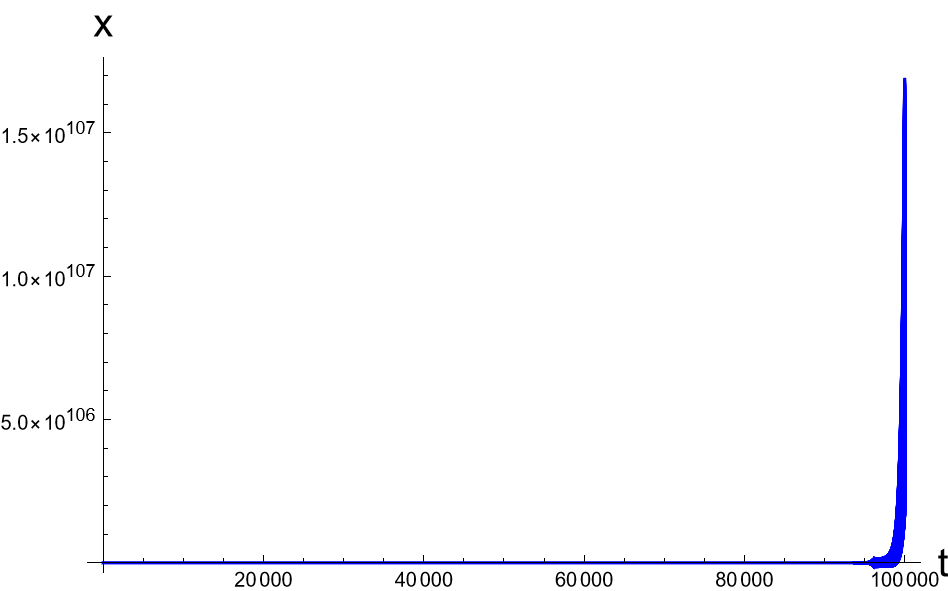} \label{fig1b}} \hspace{0.3cm}
	\caption{Solutions of system (\ref{3}) with $\alpha_1=0.2$ and $\alpha_2=0.6$ for various values of $r$ inside and outside the stable region.} \label{fig1}
\end{figure}

\begin{figure}[H]
	\centering
	\subfigure[$r=0.8$.]
	{\includegraphics[height=2.2in,keepaspectratio,width=2.2in]{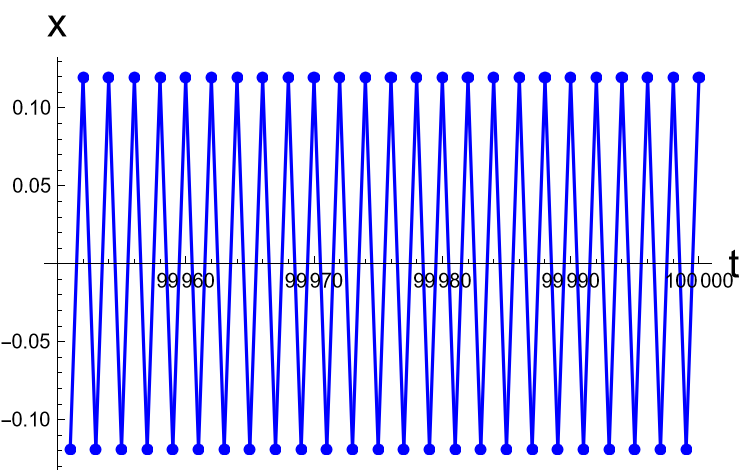}\label{fig2a}}\hspace{0.2cm}
	\subfigure[$r=-0.25$.]
	{\includegraphics[height=2.2in,keepaspectratio,width=2.2in]{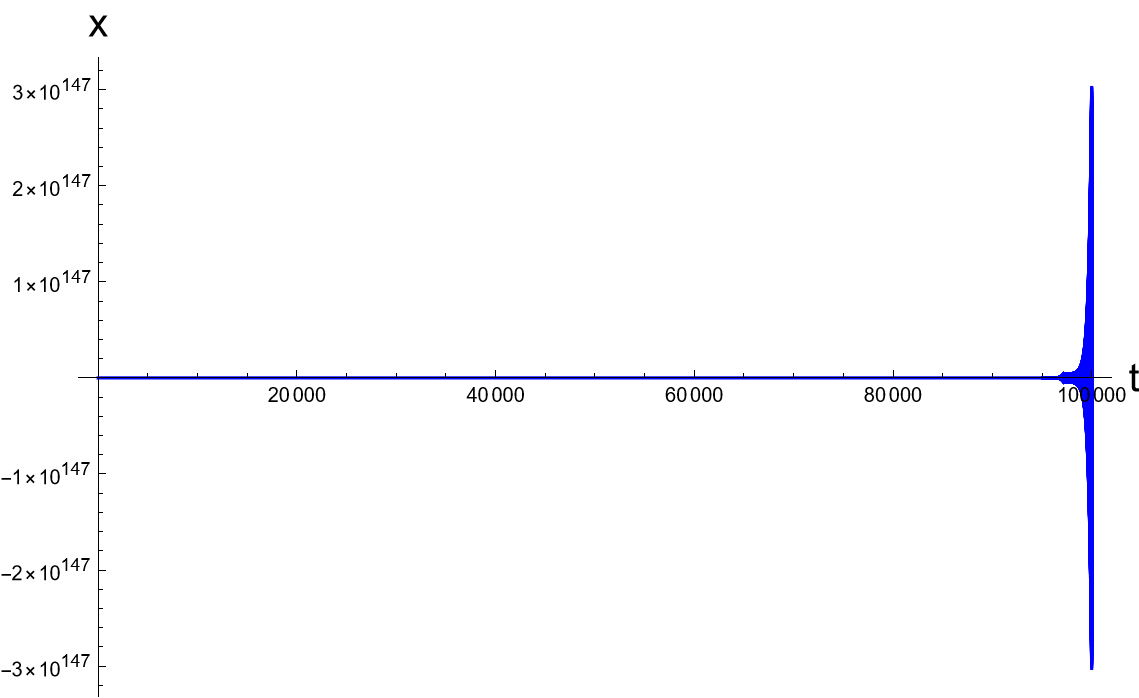}\label{fig2b}}\hspace{0.3cm}
	\caption{ Oscillations of the system (\ref{3}) with $\alpha_1=0.3$ and $\alpha_2=0.9$ for various values of $r$ giving stable and unstable oscillations.} \label{fig2}
\end{figure}

We observed that, $1-2^{\text{min}\{\alpha_1,\alpha_2\}}<r<1$ then both components $a(t)$ and $b(t)$ of the system (\ref{4}) become stable. Suppose $\lim_{t\to\infty} \frac{a(t)}{x(0)}=a_{*}$ and $\lim_{t\to\infty} \frac{b(t)}{x(0)}=b_{*}$. The system (\ref{3}) under this condition oscillates between $a_{*}$ and $b_{*}$. It becomes asymptotic periodic to $\{a_{*},b_{*}\}$. We provide few particular examples below.
\begin{ex}
	Let $\alpha_1=0.2$ and $\alpha_2=0.6$. In this case, the system (\ref{3}) with period-2 order is stable for $r \in (-0.148698,1)$. We take $r=0.4$ inside the stable region and the initial condition as $x(0)=0.1$. The Fig. \ref{fig1a} shows that the system (\ref{3}) oscillates between $-0.2057$ and $0.2060$. On the other hand, if we take $r=1.1$ outside the stable region, the system (\ref{3}) is unbounded (cf. Fig. \ref{fig1b}).
\end{ex}
\begin{ex}
	Consider $\alpha_1=0.3$ and $\alpha_2=0.9$. For $r \in (-0.231144,1)$, the system (\ref{3}) is stable in this case. Let $r=0.8$ be a sample value inside the stable region and the initial condition $x(0)=0.1$. Thus, the system (\ref{3}) is periodic and oscillates between  $-0.1192$ and $0.1192$ (cf. Fig. \ref{fig2a}). However, for $r=-0.25$ outside the stable region, the Fig. \ref{fig2b} shows that system (\ref{3}) is unbounded.
\end{ex}

\section{Maps with period-3 fractional order} \label{per3}
Consider the map (\ref{3}) with
\begin{equation*}
	\alpha(t)=
	\begin{cases}
		\alpha_1,\quad  \text{if t=3n} \\
		\alpha_2, \quad \text{if t=3n+1} \\
		\alpha_3, \quad \text{if t=3n+2}
	\end{cases}
\end{equation*}
where $0 < \alpha_1 < \alpha_2 < \alpha_3 \le 1$,  $t \in \mathbb{N} \cup \{0\}$ and $n \in \mathbb{N} \cup \{0\}$.\\
Using $a(t)=x(3t)$, $b(t)=x(3t+1)$ and $c(t)=x(3t+2)$, we can write an equivalent representation of the system (\ref{3}) as 
\begin{equation}
	\begin{split}
		a(t+1) &= x(0)+\sum_{k=0}^t \frac{\Gamma(3t+2-3k+\alpha_3)}{\Gamma(\alpha_3) \Gamma(3t-3k+3)} (r-1) a(k)\\
		& + \sum_{k=0}^t \frac{\Gamma(3t+1-3k+\alpha_3)}{\Gamma(\alpha_3) \Gamma(3t-3k+2)} (r-1) b(k)\\
		& + \sum_{k=0}^t \frac{\Gamma(3t-3k+\alpha_3)}{\Gamma(\alpha_3) \Gamma(3t-3k+1)} (r-1) c(k),\\
	b(t+1) &= x(0)+\sum_{k=0}^{t+1} \frac{\Gamma(3t+3-3k+\alpha_1)}{\Gamma(\alpha_1) \Gamma(3t-3k+4)} (r-1) a(k)\\
	& + \sum_{k=0}^t \frac{\Gamma(3t+2-3k+\alpha_1)}{\Gamma(\alpha_1) \Gamma(3t-3k+3)} (r-1) b(k)\\
	& + \sum_{k=0}^t \frac{\Gamma(3t-3k+1+\alpha_1)}{\Gamma(\alpha_1) \Gamma(3t-3k+2)} (r-1) c(k),\\
		c(t+1) &= x(0)+\sum_{k=0}^{t+1} \frac{\Gamma(3t+4-3k+\alpha_2)}{\Gamma(\alpha_2) \Gamma(3t-3k+5)} (r-1) a(k)\\
	& + \sum_{k=0}^{t+1} \frac{\Gamma(3t+3-3k+\alpha_2)}{\Gamma(\alpha_2) \Gamma(3t-3k+4)} (r-1) b(k)\\
	& + \sum_{k=0}^t \frac{\Gamma(3t-3k+2+\alpha_2)}{\Gamma(\alpha_2) \Gamma(3t-3k+3)} (r-1) c(k). \label{11}
	\end{split}
\end{equation}
If we define $\psi_{\alpha_3}(n_0)=\Tilde{\phi}_{\alpha_3}(3n_0+2)$, $\psi_{\alpha_3}^{'} (n_0)=\Tilde{\phi}_{\alpha_3}(3n_0+1)$, $\psi_{\alpha_3}^{''} (n_0)=\Tilde{\phi}_{\alpha_3}(3n_0)$, $\psi_{\alpha_1}(n_0)=\Tilde{\phi}_{\alpha_1}(3n_0+3)$, $\psi_{\alpha_1}^{'} (n_0)=\Tilde{\phi}_{\alpha_1}(3n_0+2)$, $\psi_{\alpha_1}^{''} (n_0)=\Tilde{\phi}_{\alpha_1}(3n_0+1)$, $\psi_{\alpha_2}(n_0)=\Tilde{\phi}_{\alpha_2}(3n_0+4)$, $\psi_{\alpha_2}^{'} (n_0)=\Tilde{\phi}_{\alpha_2}(3n_0+3)$ and $\psi_{\alpha_2}^{''} (n_0)=\Tilde{\phi}_{\alpha_2}(3n_0+2)$, then the system (\ref{11}) can be expressed as  
\begin{equation}
	\begin{split}
		a(t+1) &= x(0)+(r-1)\left[(\psi_{\alpha_3}\ast a)(t)\right]+(r-1)\left[(\psi_{\alpha_3}^{'}\ast b)(t)\right]\\
		&+(r-1)\left[ (\psi_{\alpha_3}^{''}\ast c)(t)\right],\\
		b(t+1) &= x(0)+ (r-1)\left[(\psi_{\alpha_1}\ast a)(t)\right]\\
		&+(r-1) a(t+1)+(r-1) \left[ (\psi_{\alpha_1}^{'}\ast b)(t)\right] +(r-1) \left[ (\psi_{\alpha_1}^{''}\ast c)(t)\right] ,\\
		c(t+1) &= x(0)+ (r-1) \left[ (\psi_{\alpha_2}\ast a)(t)\right]+(r-1) a(t+1)\\ &+(r-1) \left[ (\psi_{\alpha_2}^{'}\ast b)(t)\right] +(r-1) b(t+1)+(r-1) \left[ (\psi_{\alpha_2}^{''}\ast c)(t)\right]. \label{12}
	\end{split}
\end{equation}

The Z-transforms give
\begin{equation}
\begin{split}
Z(\psi_{\alpha_3}(t)) &= \frac{(z^{1/3})^2[(1-z^{-1/3})^{-\alpha_3}+\omega(1-\frac{\omega}{z^{1/3}})^{-\alpha_3}+\omega^2(1-\frac{\omega^2}{z^{1/3}})^{-\alpha_3}]}{3},\\
Z(\psi_{\alpha_3}^{'}(t)) &= \frac{z^{1/3}[(1-z^{-1/3})^{-\alpha_3}+\omega^2(1-\frac{\omega}{z^{1/3}})^{-\alpha_3}+\omega(1-\frac{\omega^2}{z^{1/3}})^{-\alpha_3}]}{3},\\
Z(\psi_{\alpha_3}^{''}(t)) &= \frac{[(1-z^{-1/3})^{-\alpha_3}+(1-\frac{\omega}{z^{1/3}})^{-\alpha_3}+(1-\frac{\omega^2}{z^{1/3}})^{-\alpha_3}]}{3},\\
Z(\psi_{\alpha_1}(t)) &= \frac{z[(1-z^{-1/3})^{-\alpha_1}+(1-\frac{\omega}{z^{1/3}})^{-\alpha_1}+(1-\frac{\omega^2}{z^{1/3}})^{-\alpha_1}]}{3}-z,\\
Z(\psi_{\alpha_1}^{'}(t)) &= \frac{(z^{1/3})^2[(1-z^{-1/3})^{-\alpha_1}+\omega(1-\frac{\omega}{z^{1/3}})^{-\alpha_1}+\omega^2(1-\frac{\omega^2}{z^{1/3}})^{-\alpha_1}]}{3},\\
Z(\psi_{\alpha_1}^{''}(t)) &= \frac{z^{1/3}[(1-z^{-1/3})^{-\alpha_1}+\omega^2(1-\frac{\omega}{z^{1/3}})^{-\alpha_1}+\omega(1-\frac{\omega^2}{z^{1/3}})^{-\alpha_1}]}{3},\\
Z(\psi_{\alpha_2}(t)) &= \frac{(z^{1/3})^4[(1-z^{-1/3})^{-\alpha_2}+\omega^2(1-\frac{\omega}{z^{1/3}})^{-\alpha_2}+\omega(1-\frac{\omega^2}{z^{1/3}})^{-\alpha_2}]}{3}-z \alpha_2,\\
Z(\psi_{\alpha_2}^{'}(t)) &= \frac{z[(1-z^{-1/3})^{-\alpha_2}+(1-\frac{\omega}{z^{1/3}})^{-\alpha_2}+(1-\frac{\omega^2}{z^{1/3}})^{-\alpha_2}]}{3}-z,\\
Z(\psi_{\alpha_2}^{''}(t)) &= \frac{(z^{1/3})^2[(1-z^{-1/3})^{-\alpha_2}+\omega(1-\frac{\omega}{z^{1/3}})^{-\alpha_2}+\omega^2(1-\frac{\omega^2}{z^{1/3}})^{-\alpha_2}]}{3},\\
Z(a(t+1)) &= zA(z)-za(0), \quad Z(b(t+1))=zB(z)-zb(0), \quad \text{and} \\ 
 Z(c(t+1)) &= zC(z)-zc(0), 
\end{split}
 \label{13}
\end{equation}
where $Z(a(t))=A(z)$, $Z(b(t))=B(z)$, $Z(c(t))=C(z)$, $a(0)=x(0)$, $b(0)=x(1)$, $c(0)=x(2)$.\\
Using (\ref{13}) and applying Z-transform to the system (\ref{12}), we obtain\\
\begin{equation*}
	\begin{split}
		&((r-1)\frac{(z^{1/3})^2[(1-z^{-1/3})^{-\alpha_3}+\omega(1-\frac{\omega}{z^{1/3}})^{-\alpha_3}+\omega^2(1-\frac{\omega^2}{z^{1/3}})^{-\alpha_3}]}{3}-z)A(z)\\
		&+((r-1)\frac{z^{1/3}[(1-z^{-1/3})^{-\alpha_3}+\omega^2(1-\frac{\omega}{z^{1/3}})^{-\alpha_3}+\omega(1-\frac{\omega^2}{z^{1/3}})^{-\alpha_3}]}{3})B(z)\\
		&+((r-1)\frac{[(1-z^{-1/3})^{-\alpha_3}+(1-\frac{\omega}{z^{1/3}})^{-\alpha_3}+(1-\frac{\omega^2}{z^{1/3}})^{-\alpha_3}]}{3})C(z)=R_1, 
	\end{split}
\end{equation*}
\begin{equation*}
	\begin{split}
		&((r-1)\frac{z[(1-z^{-1/3})^{-\alpha_1}+(1-\frac{\omega}{z^{1/3}})^{-\alpha_1}+(1-\frac{\omega^2}{z^{1/3}})^{-\alpha_1}]}{3})A(z)\\
		&+((r-1)\frac{(z^{1/3})^2[(1-z^{-1/3})^{-\alpha_1}+\omega(1-\frac{\omega}{z^{1/3}})^{-\alpha_1}+\omega^2(1-\frac{\omega^2}{z^{1/3}})^{-\alpha_1}]}{3}-z)B(z)\\
		&+((r-1)\frac{z^{1/3}[(1-z^{-1/3})^{-\alpha_1}+\omega^2(1-\frac{\omega}{z^{1/3}})^{-\alpha_1}+\omega(1-\frac{\omega^2}{z^{1/3}})^{-\alpha_1}]}{3})C(z)=R_2, 
	\end{split}
\end{equation*}
\begin{equation}
	\begin{split}
		&((r-1)\frac{(z^{1/3})^4[(1-z^{-1/3})^{-\alpha_2}+\omega^2(1-\frac{\omega}{z^{1/3}})^{-\alpha_2}+\omega(1-\frac{\omega^2}{z^{1/3}})^{-\alpha_2}]}{3})A(z)\\
		&+((r-1)\frac{z[(1-z^{-1/3})^{-\alpha_2}+(1-\frac{\omega}{z^{1/3}})^{-\alpha_2}+(1-\frac{\omega^2}{z^{1/3}})^{-\alpha_2}]}{3})B(z)\\
		&+((r-1)\frac{(z^{1/3})^2[(1-z^{-1/3})^{-\alpha_2}+\omega(1-\frac{\omega}{z^{1/3}})^{-\alpha_2}+\omega^2(1-\frac{\omega^2}{z^{1/3}})^{-\alpha_2}]}{3}-z)C(z)\\
        &=R_3. \label{14}
	\end{split}
\end{equation}	
where $R_1=-\frac{x(0)}{1-z^{-1}}-zx(0)$, $R_2=(r-1)zx(0)-\frac{x(0)}{1-z^{-1}}-zx(1)$ and $R_3=(r-1)z\alpha_2 x(0)-\frac{x(0)}{1-z^{-1}}+(r-1)z x(1)-z x(2)$ are remainder terms containing initial values and $\omega= e^{\frac{2 \pi i}{3}}$.
\subsection{The Characteristic Equation and Stable Region} \label{}
The characteristic equation of the system (\ref{11}) can be obtained by equating the determinant of coefficients of the terms A(z), B(z) and C(z) in the system (\ref{14}) to zero. We provided the expression for the characteristic equation in the Data set-1 available at the link \url{https://github.com/Janardhan3233/VOM.git}. 
\par As in Section \ref{3.1}, by substituting $z=\pm 1$ in the characteristic equation of the system (\ref{11}), the boundary of the system's stable region can be found. Substituting $z=z^{1/3}=-1$ in the characteristic equation of the system (\ref{11}), we get the boundaries of stable region as below :
\begin{equation}
	\left( \frac{r-1}{3}\right)^3
	\left[A_1(B_2 C_3-C_2 B_3)-B_1(A_2 C_3-C_2 A_3)+C_1(A_2 B_3-B_2 A_3)\right]=0. \label{15}
\end{equation}
where $A_1=2^{-\alpha_3}+\frac{3}{-1+r}-\cos(\frac{\pi \alpha_3}{3})+\sqrt{3} \sin(\frac{\pi \alpha_3}{3})$, $B_1=-2^{-\alpha_3}+\cos(\frac{\pi \alpha_3}{3})+\sqrt{3} \sin(\frac{\pi \alpha_3}{3})$, $C_1=2^{-\alpha_3}+2 \cos(\frac{\pi \alpha_3}{3})$, $A_2=-2^{-\alpha_1}-2 \cos(\frac{\pi \alpha_1}{3})$, $B_2=2^{-\alpha_1}+\frac{3}{-1+r}-\cos(\frac{\pi \alpha_1}{3})+\sqrt{3} \sin(\frac{\pi \alpha_1}{3})$, $C_2=-2^{-\alpha_1}+\cos(\frac{\pi \alpha_1}{3})+\sqrt{3} \sin(\frac{\pi \alpha_1}{3})$, $A_3=2^{-\alpha_2}-\cos(\frac{\pi \alpha_2}{3})-\sqrt{3} \sin(\frac{\pi \alpha_2}{3})$, $B_3=-2^{-\alpha_2}-2 \cos(\frac{\pi \alpha_2}{3})$, and $C_3=2^{-\alpha_2}+\frac{3}{-1+r}-\cos(\frac{\pi \alpha_2}{3})+\sqrt{3} \sin(\frac{\pi \alpha_2}{3})$.\\
Solving the equation (\ref{15}), we get $r=1$ as the right boundary and another expression of $r$ given in Data set-2 (available at the above link) as the left boundary of the stable region of the system (\ref{11}).
\par
\begin{figure}[H]
	\includegraphics[height=4.2in,keepaspectratio,width=4.2in]{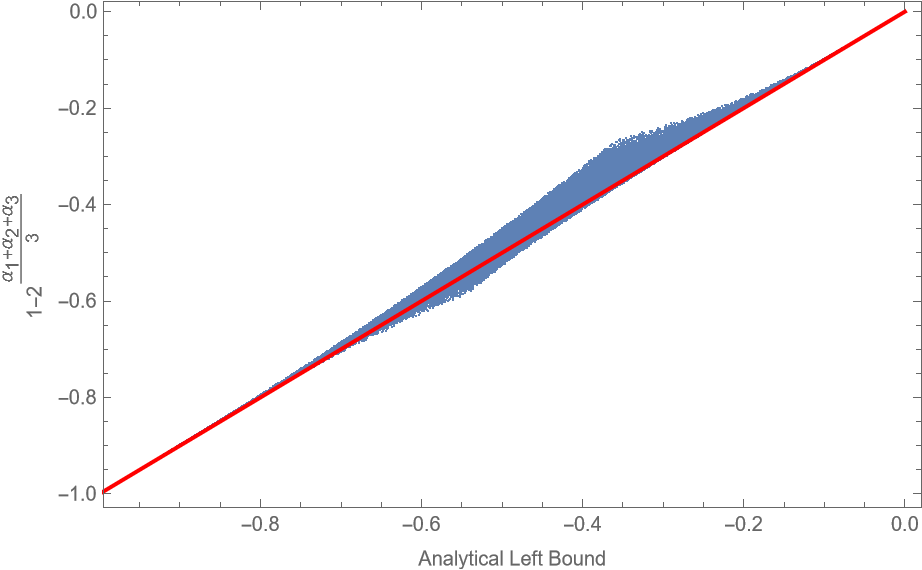}
	\caption{Analytical value left bound versus $1-2^{\frac{\alpha_1+\alpha_2+\alpha_3}{3}}$ for 1 lakh random triples $(\alpha_1,\alpha_2,\alpha_3)$ with $\alpha_j\in (0,1)$, $j=1,2,3$.}
	\label{fig3}
	\end{figure}

 In Fig. \ref{fig3}, we sketch the left boundary of the stable region versus $1-2^{\frac{\alpha_1+\alpha_2+\alpha_3}{3}}$. It is observed that the value $1-2^{\frac{\alpha_1+\alpha_2+\alpha_3}{3}}$ is the best approximation to the left boundary.

\subsection{Validation of analytic results  with Numerical Experiments}
\begin{figure}[H]
	\centering
	\subfigure[$r=0.8$.]
	{\includegraphics[height=2.2in,keepaspectratio,width=2.2in]{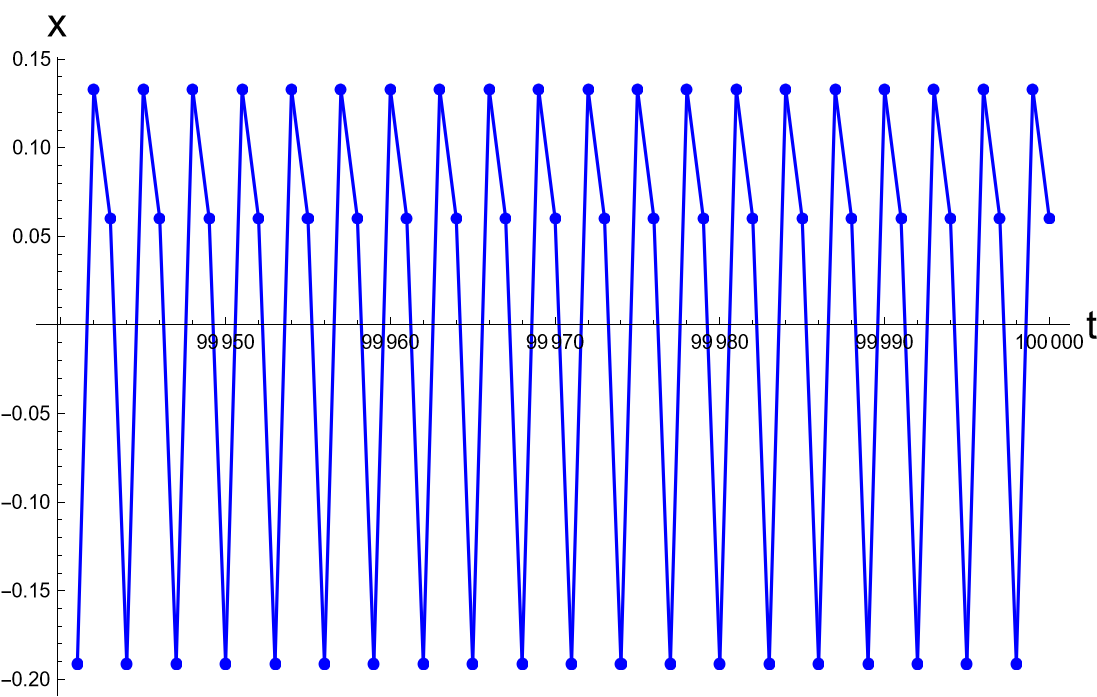}\label{fig4a}} \hspace{0.2cm}
	\subfigure[$r=1.1$.]
	{\includegraphics[height=2.2in,keepaspectratio,width=2.2in]{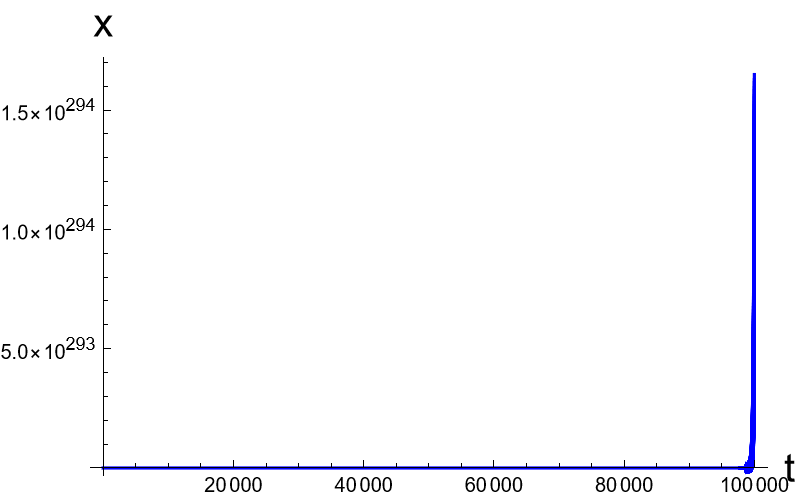} \label{fig4b}} \hspace{0.3cm}
	\caption{Behavior of the system (\ref{3}) with $\alpha_1=0.2$, $\alpha_2=0.4$ and $\alpha_3=0.6$ for different values of $r$ inside and outside the stable region.} \label{fig4}
\end{figure}

\begin{figure}[H]
	\centering
	\subfigure[$r=-0.1$.]
	{\includegraphics[height=2.2in,keepaspectratio,width=2.2in]{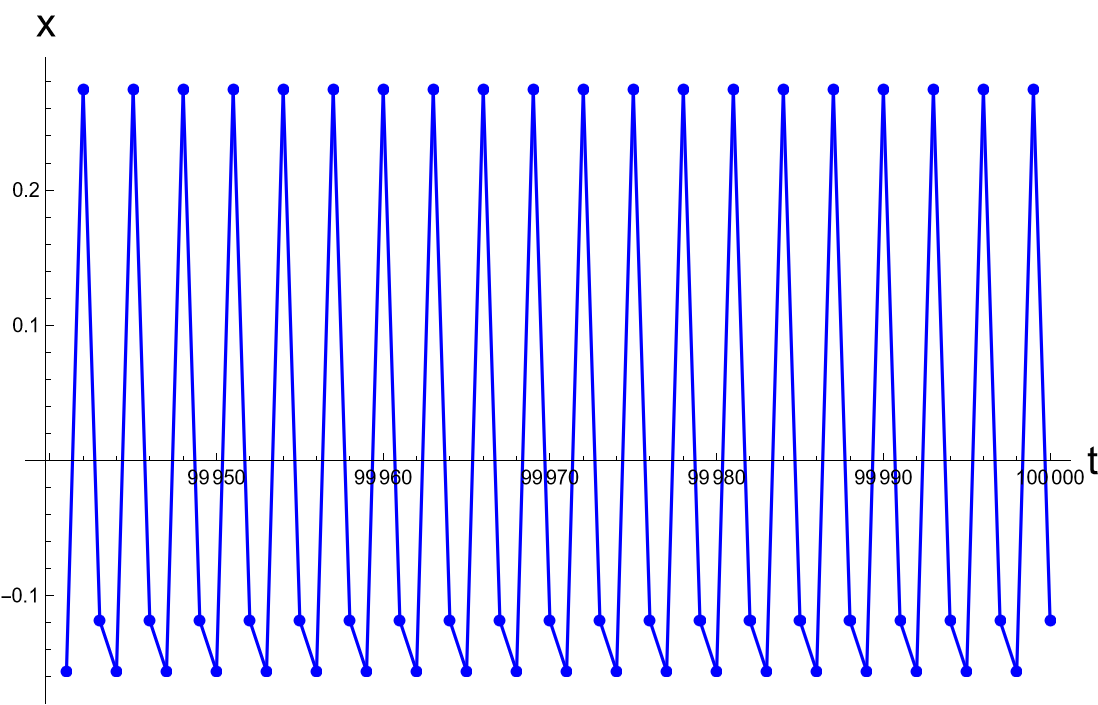}\label{fig5a}}\hspace{0.2cm}
	\subfigure[$r=-0.5$.]
	{\includegraphics[height=2.2in,keepaspectratio,width=2.2in]{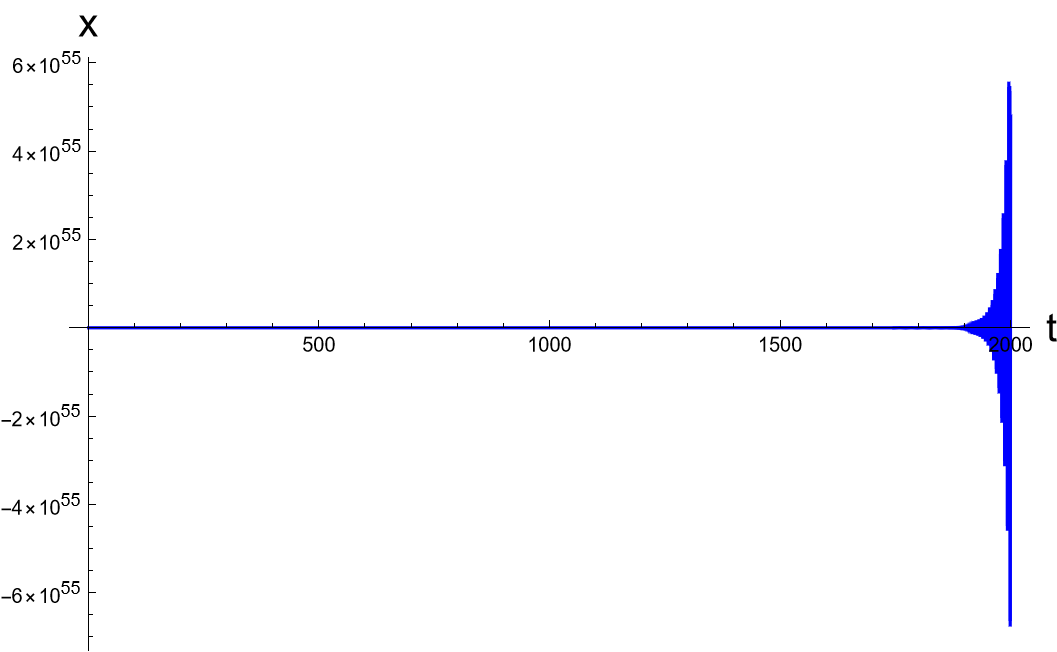}\label{fig5b}}\hspace{0.3cm}
	\caption{Solutions of system (\ref{3}) with $\alpha_1=0.1$, $\alpha_2=0.5$ and $\alpha_3=0.9$ for different values of $r$.} \label{fig5}
\end{figure}

 The system (\ref{3}) exhibits asymptotic period-3 oscillations, if we take the $r$ values in the stable region. Let $\lim_{t\to\infty} \frac{a(t)}{x(0)}=a_{*}$, $\lim_{t\to\infty} \frac{b(t)}{x(0)}=b_{*}$ and $\lim_{t\to\infty} \frac{c(t)}{x(0)}=c_{*}$, for some numbers $a_{*}$, $b_{*}$, and $c_{*}$. A few specific examples are given below.
\begin{ex}
	Assume that $\alpha_1=0.2$, $\alpha_2=0.4$ and $\alpha_3=0.6$. For $r \in (-0.322259,1)$, the system (\ref{11}) is stable in this case. Consider $r=0.8$ inside the stable region and the initial condition $x(0)=0.1$. According to Fig. \ref{fig4a}, the system (\ref{3}) is asymptotically periodic to $\{ -0.1913, 0.1331, 0.0602 \}$. On the other hand, the system (\ref{3}) is unbounded if $r=1.1$, outside the stable region (see Fig. \ref{fig4b}).
\end{ex}
\begin{ex}
	Consider $\alpha_1=0.1$, $\alpha_2=0.5$ and $\alpha_3=0.9$. For $r=-0.1$ inside the stable region $(-0.432383,1)$ and the initial condition $x(0)=0.1$, the system (\ref{3}) exhibits asymptotically periodic-3 oscillations between $-0.1186$, $-0.1558$, and $0.2744$ (cf. Fig. \ref{fig5a}). However, for $r=-0.5$ outside the stable region, the system (\ref{3}) is unbounded (cf. Fig. \ref{fig5b}).
\end{ex}

\section{Maps with higher periodic orders}\label{pern}
\begin{figure}[H]
	\includegraphics[height=4.2in,keepaspectratio,width=4.2in]{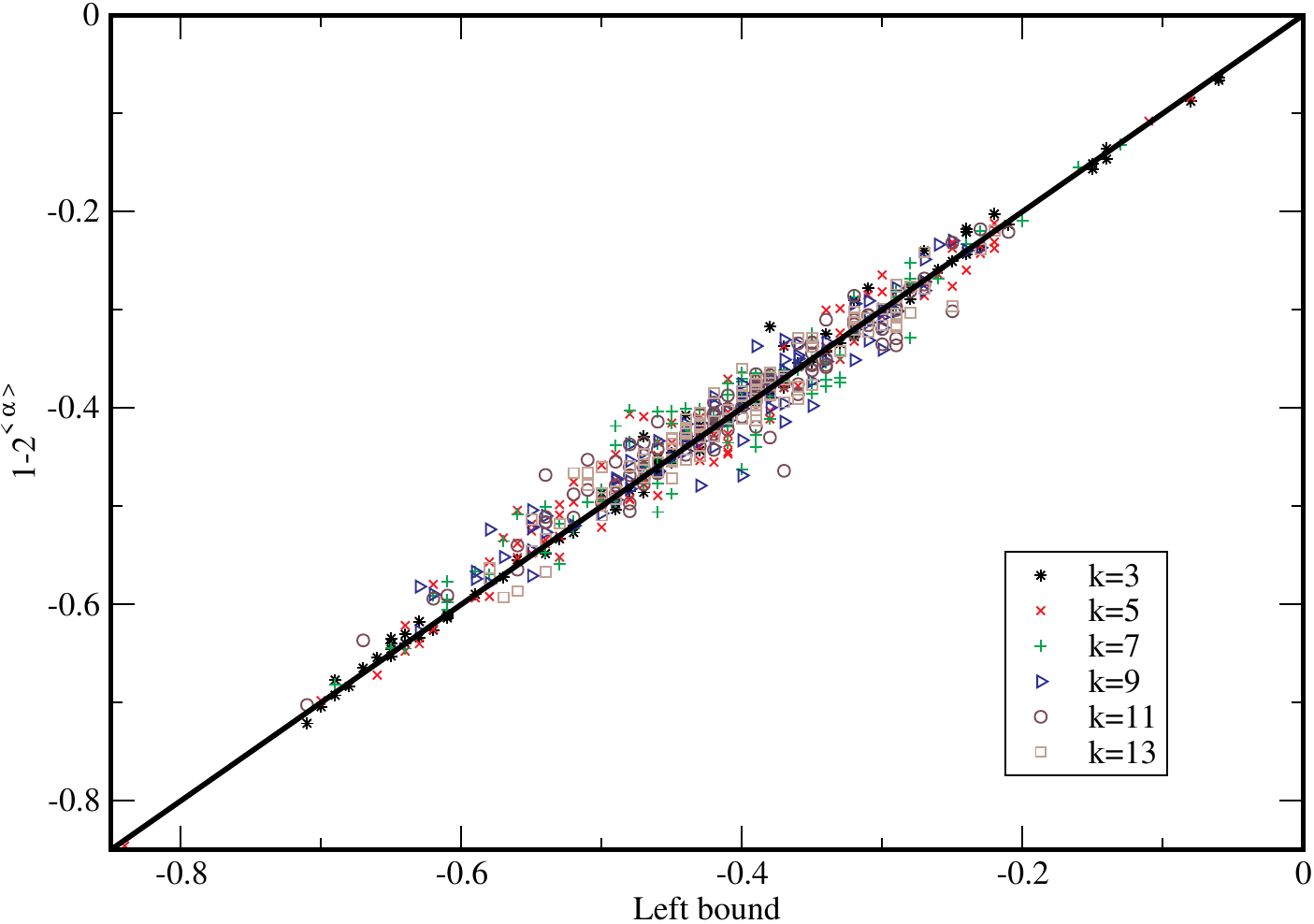}
	\caption{Experimental value of left boundary of stable region after $3.8\times 10^8$ time-steps for 100 different random configurations versus $1-2^{\langle\alpha\rangle}$ for the periods 3 (Black), 5 (Red), 7 (Green), 9 (Blue), 11 (Magenta) and 13 (Brown). }
	\label{fig6}
\end{figure}

  Consider a case of Eq. (\ref{3}) where $\alpha(t)$ is a
periodic function with period T. 
\begin{equation*}
	\alpha(t)=\alpha_k \text{ where  } k= \mod(t,T)+1,
\end{equation*}
where  $t \in \mathbb{N} \cup \{0\}$ and $T\in \mathbb{N}$.\\
Let $\langle \alpha \rangle $ be the mean value of $\alpha(t)$,
\it{ i.e.} $\langle \alpha  \rangle ={\frac{1}{T}} \sum_{i=1}^T \alpha_i $. We propose that this mean value is a good indicator of the stability of this system for $T\ne 2, 4, 6$ and 8. For $T=2$, the $\alpha_{min}=\min(\alpha_1,\alpha_2)$  dictates the stability of the system. The system is stable if $1-2^{\alpha_{min}}<r<1$. For $T=3$, the theoretical left stability bound is very close to $1-2^{\langle\alpha\rangle}$ and the right bound is 1. For odd values of $T$, the stability is well approximated by the range $1-2^{\langle\alpha\rangle} <r< 1$ (see Fig. \ref{fig6}).
    
\begin{figure}[H]
    \centering
    \includegraphics[width=1\textwidth]{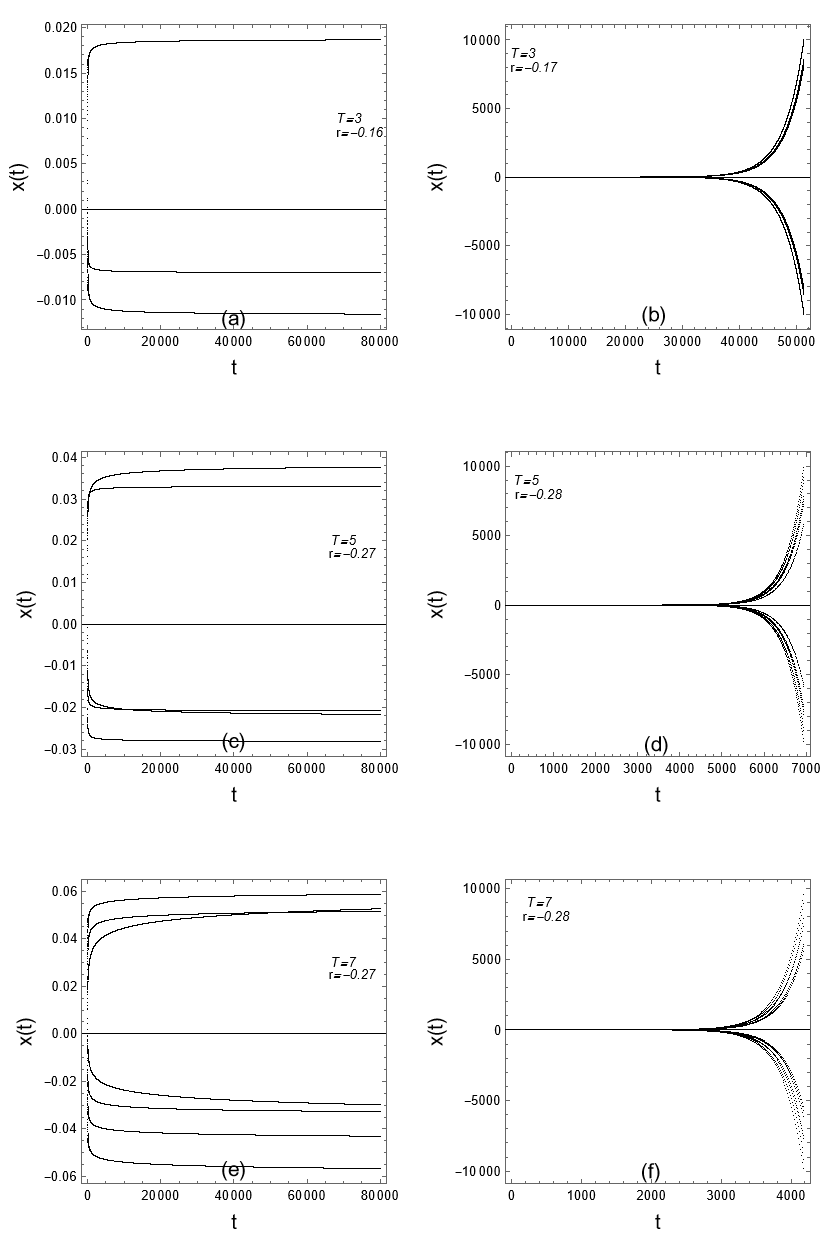}
    \caption{ The time series is plotted for $T=3$ ($1-2^{\langle\alpha\rangle}=-0.1638$) and a)$r=-0.16$ and b)$r=-0.17$.  The time series is plotted for $T=5$ ($1-2^{\langle\alpha\rangle}= -0.2846$) for c) r=-0.28 and d) r=-0.27.
    The time series is plotted for $T=7$ ($1-2^{\langle\alpha\rangle}=-0.2987 $) for e) $r= -0.27$ f)$r=-0.28$. The $\alpha$ values are chosen randomly.}  \label{fig7}
\end{figure}

\begin{ex}
\par For $T=3$, Fig. \ref{fig7}a shows that the system is in stable (periodic) orbit for $r=-0.16$. The system is unstable for $r=-0.17$ outside the stable region (see Fig. \ref{fig7}b). In this  case, The threshold value is near the expected value $r=1-2^{\langle\alpha\rangle}=-0.1638$. 
\par For $T=5$, the system is in stable (periodic) orbit for $r=-0.27$, as seen in Fig. \ref{fig7}c. For $r=-0.28$ beyond the stable region, the system is unstable (see Fig. \ref{fig7}d). The threshold value in this case is higher than the predicted value $r=-0.2846$.
\par For $T=7$, Fig. \ref{fig7}e indicates that the system is in stable (periodic) orbit for $r=-0.27$. Outside of the stable region (see Fig. \ref{fig7}f), the system is unstable for $r=-0.28$, where the threshold value is higher than the expected value $r=-0.2987$. The $\alpha$ values are chosen randomly and the initial condition as $x(0)=0.01$ in all above cases. 
We observe that the numerical value of threshold is within
$0.02$ of the expected value for $T=3, 5$ and $7$.
\end{ex}

\begin{figure}[H]
    \centering
    \includegraphics[width=1\textwidth]{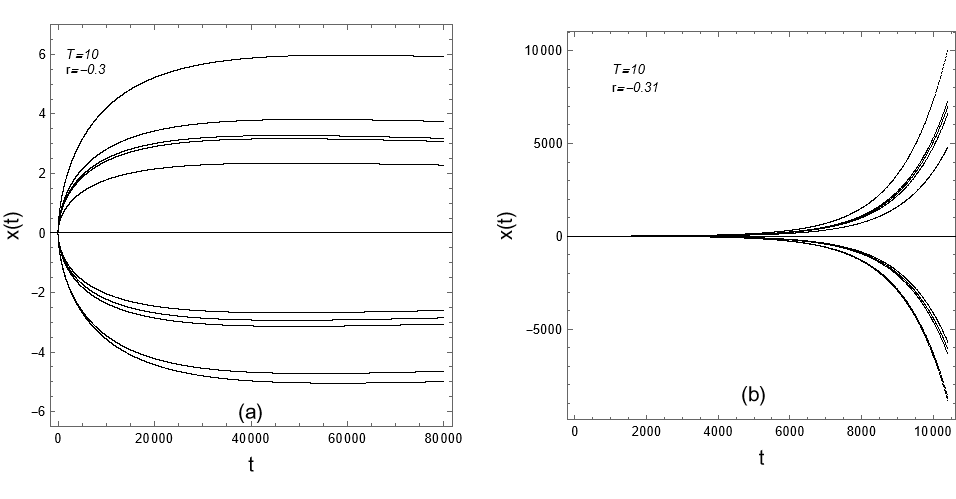} 
    \caption{For $T=10$, the time series is plotted for variable order map with 10 period for a)r=-0.31 and b)r=-0.3. The $\alpha$ values are chosen randomly. The expected value is -0.4052. Thus the threshold differs considerably from the expected value.} \label{fig8}
\end{figure}

\begin{ex}
 Consider $T=10$. The $\alpha$ values are chosen randomly and the initial condition as $x(0)=0.01$. According to fig. \ref{fig8}a, for $r=-0.30$, the system is in a stable (periodic) orbit. The system is unstable for $r=-0.31$ outside the stable region (see Fig. \ref{fig8}b). The threshold value in this instance deviates significantly from the expected value $r=-0.4052$.
\end{ex}

\section{Conclusion} \label{con.}
	In this work, we have studied the variable order fractional maps with periodic variational order. The cases of periods 2 and 3 have been tackled analytically and explicit stability conditions are derived for a given sequence. The upper bound is unity while the lower bound depends on the sequence. For period-2, we show that the stability condition
    is dictated by the smaller of the two fractional orders.
    For period-3, the stability conditions are complicated. But our numerical experiments show that the analytic value of the lower bound is close to that for 
    fractional order map of constant order $\langle \alpha\rangle $ where $\langle \alpha\rangle$ is the mean 
    value of the fractional orders. The analytic bound is corroborated by numerical experiments in all the cases. 

    For even periods larger than 2, it is difficult to get an analytic estimate. The numerical estimate of the lower bound does not show a clear pattern. A naive guess that at least for large values of even $T$, we should get results similar to odd values, and the stability condition should be dictated by $\langle\alpha\rangle$ does not seem to hold. Why such a simplistic guess works for odd period  is an open question.
    Another open question is why this condition does not hold for even periods.
	
	\section*{Acknowledgments}
	P. M. Gade thanks DST-SERB for financial assistance (Ref. CRG/2020/003993). Ch. Janardhan thanks University Grants Commission, New Delhi, India for financial support (No. F.14-34/2011(CPP-II)).
	
	\section*{Datasets}
	All the datasets referenced in the context of this research work are stored at the link \url{https://github.com/Janardhan3233/VOM.git}
	
	%
	
	

\bibliographystyle{plain}      
\bibliography{ref.bib}

\begin{thebibliography}{10}

\bibitem{abu2013asymptotic}
Raghib Abu-Saris and Qasem Al-Mdallal.
\newblock On the asymptotic stability of linear system of fractional-order
  difference equations.
\newblock {\em Fractional Calculus and Applied Analysis}, 16(3):613--629, 2013.

\bibitem{albasheir2023theoretical}
Nafisa~A Albasheir, Ammar Alsinai, Azmat Ullah~Khan Niazi, Ramsha Shafqat,
  Romana, Mohammed Alhagyan, and Ameni Gargouri.
\newblock A theoretical investigation of caputo variable order fractional
  differential equations: existence, uniqueness, and stability analysis.
\newblock {\em Computational and Applied Mathematics}, 42(8):367, 2023.

\bibitem{almatroud2024new}
Othman~Abdullah Almatroud, Amina-Aicha Khennaoui, Adel Ouannas, Saleh
  Alshammari, and Sahar Albosaily.
\newblock A new fractional discrete memristive map with variable order and
  hidden dynamics.
\newblock {\em Fractal and Fractional}, 8(6):322, 2024.

\bibitem{bastos2011discrete}
Nuno~RO Bastos, Rui~AC Ferreira, and Delfim~FM Torres.
\newblock Discrete-time fractional variational problems.
\newblock {\em Signal Processing}, 91(3):513--524, 2011.

\bibitem{bensid2022discrete}
Souad Bensid~Ahmed, Adel Ouannas, Mohammed Al~Horani, and Giuseppe Grassi.
\newblock The discrete fractional variable-order tinkerbell map: chaos, 0--1
  test, and entropy.
\newblock {\em Mathematics}, 10(17):3173, 2022.

\bibitem{bhalekar2022stability}
Sachin Bhalekar and Prashant~M Gade.
\newblock Stability analysis of fixed point of fractional-order coupled map
  lattices.
\newblock {\em Communications in Nonlinear Science and Numerical Simulation},
  113:106587, 2022.

\bibitem{cao2016high}
Jianxiong Cao and Yanan Qiu.
\newblock A high order numerical scheme for variable order fractional ordinary
  differential equation.
\newblock {\em Applied Mathematics Letters}, 61:88--94, 2016.

\bibitem{vcermak2015explicit}
Jan {\v{C}}erm{\'a}k, Istv{\'a}n Gy{\H{o}}ri, and Lud{\u{e}}k Nechv{\'a}tal.
\newblock On explicit stability conditions for a linear fractional difference
  system.
\newblock {\em Fractional Calculus and Applied Analysis}, 18:651--672, 2015.

\bibitem{chen2014numerical}
Shiping Chen, Fawang Liu, and Kevin Burrage.
\newblock Numerical simulation of a new two-dimensional variable-order
  fractional percolation equation in non-homogeneous porous media.
\newblock {\em Computers \& Mathematics with Applications}, 68(12):2133--2141,
  2014.

\bibitem{coimbra2003mechanics}
Carlos~FM Coimbra.
\newblock Mechanics with variable-order differential operators.
\newblock {\em Annalen der Physik}, 515(11-12):692--703, 2003.

\bibitem{derakhshan2022existence}
Mohammad~Hossein Derakhshan.
\newblock Existence, uniqueness, ulam--hyers stability and numerical simulation
  of solutions for variable order fractional differential equations in fluid
  mechanics.
\newblock {\em Journal of Applied Mathematics and Computing}, 68(1):403--429,
  2022.

\bibitem{elaydi2005systems}
Saber Elaydi.
\newblock Systems of linear difference equations.
\newblock {\em An Introduction to Difference equations}, pages 117--172, 2005.

\bibitem{ferreira2011fractional}
Rui~AC Ferreira and Delfim~FM Torres.
\newblock Fractional h-difference equations arising from the calculus of
  variations.
\newblock {\em Applicable Analysis and Discrete Mathematics}, pages 110--121,
  2011.

\bibitem{fulai2011existence}
Chen Fulai, Luo Xiannan, and Zhou Yong.
\newblock Existence results for nonlinear fractional difference equation.
\newblock {\em Advances in Difference Equations}, 12, 2011.

\bibitem{hamadneh2023finite}
Tareq Hamadneh, Amel Hioual, Omar Alsayyed, Yazan~Alaya Al-Khassawneh, Abdallah
  Al-Husban, and Adel Ouannas.
\newblock Finite time stability results for neural networks described by
  variable-order fractional difference equations.
\newblock {\em Fractal and Fractional}, 7(8):616, 2023.

\bibitem{hirsch2012differential}
Morris~W Hirsch, Stephen Smale, and Robert~L Devaney.
\newblock {\em Differential equations, dynamical systems, and an introduction
  to chaos}.
\newblock Academic press, 2012.

\bibitem{khan2023system}
Hasib Khan, Jehad Alzabut, Haseena Gulzar, Osman Tun{\c{c}}, and Sandra
  Pinelas.
\newblock On system of variable order nonlinear p-laplacian fractional
  differential equations with biological application.
\newblock {\em Mathematics}, 11(8):1913, 2023.

\bibitem{lorenzo2002variable}
Carl~F Lorenzo and Tom~T Hartley.
\newblock Variable order and distributed order fractional operators.
\newblock {\em Nonlinear dynamics}, 29:57--98, 2002.

\bibitem{mozyrska2019stability}
Dorota Mozyrska, Piotr Oziablo, and Ma{\l}gorzata Wyrwas.
\newblock Stability of fractional variable order difference systems.
\newblock {\em Fractional Calculus and Applied Analysis}, 22(3):807--824, 2019.

\bibitem{mozyrska2015transform}
Dorota Mozyrska, Ma{\l}gorzata Wyrwas, et~al.
\newblock The-transform method and delta type fractional difference operators.
\newblock {\em Discrete Dynamics in Nature and Society}, 2015, 2015.

\bibitem{patnaik2020applications}
Sansit Patnaik, John~P Hollkamp, and Fabio Semperlotti.
\newblock Applications of variable-order fractional operators: a review.
\newblock {\em Proceedings of the Royal Society A}, 476(2234):20190498, 2020.

\bibitem{samko1993integration}
Stefan~G Samko and Bertram Ross.
\newblock Integration and differentiation to a variable fractional order.
\newblock {\em Integral transforms and special functions}, 1(4):277--300, 1993.

\bibitem{sarwar2022existence}
Shahzad Sarwar.
\newblock On the existence and stability of variable order caputo type
  fractional differential equations.
\newblock {\em Fractal and Fractional}, 6(2):51, 2022.

\bibitem{shen2012numerical}
Shujun Shen, Fawang Liu, Jing Chen, Ian Turner, and Vo~Anh.
\newblock Numerical techniques for the variable order time fractional diffusion
  equation.
\newblock {\em Applied Mathematics and Computation}, 218(22):10861--10870,
  2012.

\bibitem{solis2018novel}
JE~Sol{\'\i}s-P{\'e}rez, JF~G{\'o}mez-Aguilar, and A~Atangana.
\newblock Novel numerical method for solving variable-order fractional
  differential equations with power, exponential and mittag-leffler laws.
\newblock {\em Chaos, Solitons \& Fractals}, 114:175--185, 2018.

\bibitem{soradi2020king}
Samaneh Soradi-Zeid, Hadi Jahanshahi, Amin Yousefpour, and Stelios Bekiros.
\newblock King algorithm: A novel optimization approach based on variable-order
  fractional calculus with application in chaotic financial systems.
\newblock {\em Chaos, Solitons \& Fractals}, 132:109569, 2020.

\bibitem{stanislawski2013stability}
R~Stanis{\l}awski and Krzysztof~J Latawiec.
\newblock Stability analysis for discrete-time fractional-order lti state-space
  systems. part i: New necessary and sufficient conditions for the asymptotic
  stability.
\newblock {\em Bulletin of the Polish Academy of Sciences. Technical Sciences},
  61(2):353--361, 2013.

\bibitem{sun2019review}
HongGuang Sun, Ailian Chang, Yong Zhang, and Wen Chen.
\newblock A review on variable-order fractional differential equations:
  mathematical foundations, physical models, numerical methods and
  applications.
\newblock {\em Fractional Calculus and Applied Analysis}, 22(1):27--59, 2019.

\bibitem{wang2019analysis}
Hong Wang and Xiangcheng Zheng.
\newblock Analysis and numerical solution of a nonlinear variable-order
  fractional differential equation.
\newblock {\em Advances in Computational Mathematics}, 45:2647--2675, 2019.

\bibitem{wu2019new}
Guo-Cheng Wu, Zhen-Guo Deng, Dumitru Baleanu, and De-Qiang Zeng.
\newblock New variable-order fractional chaotic systems for fast image
  encryption.
\newblock {\em Chaos: An Interdisciplinary Journal of Nonlinear Science},
  29(8), 2019.

\bibitem{wu2022fractional}
Guo-Cheng Wu, Chuan-Yun Gu, Lan-Lan Huang, and Dumitru Baleanu.
\newblock Fractional differential equations of variable order: existence
  results, numerical method and asymptotic stability conditions.
\newblock {\em Miskolc Mathematical Notes}, 23:485--493, 2022.

\bibitem{zayernouri2015fractional}
Mohsen Zayernouri and George~Em Karniadakis.
\newblock Fractional spectral collocation methods for linear and nonlinear
  variable order fpdes.
\newblock {\em Journal of Computational Physics}, 293:312--338, 2015.

\bibitem{zeng2015generalized}
Fanhai Zeng, Zhongqiang Zhang, and George~Em Karniadakis.
\newblock A generalized spectral collocation method with tunable accuracy for
  variable-order fractional differential equations.
\newblock {\em SIAM Journal on Scientific Computing}, 37(6):A2710--A2732, 2015.

\bibitem{zhao2015second}
Xuan Zhao, Zhi-zhong Sun, and George~Em Karniadakis.
\newblock Second-order approximations for variable order fractional
  derivatives: algorithms and applications.
\newblock {\em Journal of Computational Physics}, 293:184--200, 2015.

\end{thebibliography}

\end{document}